**Multilevel Facility Location Optimization: A Novel Integer Programming Formulation and Approaches to Heuristic Solutions**

Bahram Alidaee
School of Business Administration, The University of Mississippi, University, MS 38677, USA

balidaee@bus.olemiss.edu
662-715-1614

Haibo Wang*
A.R. Sánchez Jr. School of Business, Texas A&M International University, Laredo, TX 78045, USA

hwang@tamiu.edu
956-326-2503

**Declaration of interest statement:** The authors confirm that no known conflicts of interest are associated with this publication.

**Declaration of funding: No** funding was received.

*Corresponding author

**Multilevel Facility Location Optimization: A Novel Integer Programming Formulation and Approaches to Heuristic Solutions**

**Abstract**

We attack the 4-level facility location problem (4L-FLP), a critical component in supply chains. Foundational tasks here involve selecting markets, plants, warehouses, and distribution centers to maximize profits while considering related constraints. Based on a variation of the quadratic assignment problem, we propose a novel integer programming formula that significantly reduces the variables. Our model incorporates several realistic features, including transportation costs and upper bounds on facilities at each level. It accounts for one-time fixed costs associated with selecting each facility. To solve this complex problem, we develop and experimentally test two solution procedures: a multi-start greedy heuristic and a multi-start tabu search. We conduct extensive sensitivity analyses on the results to assess the reliability of proposed algorithms. This study contributes to improved solution methods for large-scale 4L-FLPs, providing a valuable tool for supply chain maturity.

## 1. Introduction

The focus here is on location priority in 4-level facilities (4L-FLP), a matter of far-reaching applications in various settings, including supply chain strategic structuring, transportation planning, supplier selection, manufacturing facilities, warehouses, distribution centers, and retail stores. It also has applications in many service settings, such as health care design, disaster response, telecommunication systems, postal delivery, education systems, and solid waste management. Thorough preparation cannot overlook air freight, passenger travel, forestry, oil and gas field development, and last-mile-delivery, as evident from extensive research in this field (Klose & Drexl 2005; Şahin & Süral 2007; Melo 2009, 2012; Farahani 2014; Melo. 2014; Khalifehzadeh 2015; Fattahi 2017; Mohammadi 2017; Ortiz-Astorquiza 2018; Rafiei 2018; Cortinhal 2019; Janjevic 2019; Kumar 2020; Belieres 2021; Vishnu 2021; Cao 2022; Saldanha-da-Gama 2022; Majumdar 2023; Jahani 2024; Kumar & Kumar 2024). As highlighted by Kang (2021), the practical design of delivery systems has become increasingly vital in real-world scenarios, making it a strategic priority to develop adaptive and efficient delivery networks.

The 4L-FLP involves serving a set of potential markets, also referred to as *retail stores*, *customer zones*, or *demand zones*, through a network of facilities. Each market has an associated benefit, and the objective is to maximize the total benefit by selecting the optimal facilities. The network consists of potential *plants*, *warehouses*, *distribution centers*, and *markets*. Each market is served by a *plant*, connected to a set of *warehouses,* each connected to a set of *distribution centers*. Each selected market is served by a *distribution center*. This complex network structure can be seen in a directed layered graph, illustrated in Figure 1. Graph *G* comprises a set of nodes *V*, including potential *plants P, warehouses W, distribution centers D*, and *retail stores S*, connected by a set of directed arcs *A*. Products flow from *plants* to *warehouses*: *warehouses* to *distribution centers,* and *distribution centers* to *stores* via the set of arcs *A*. While literature often uses alternative terms such as *suppliers*, *plants*, *distribution centers*, *customers*, or *agents* at each layer, we adopt the definition shown in Figure 1, which is more suitable for supply chain design, especially when entering new markets. A feasible solution to serve a selected *retail store m* involves selecting one *plant n*, one *warehouse k*, and one *distribution center j*, forming a path (*n*, *k*, *j*, *m*) in the directed graph, as depicted in Figure 1.

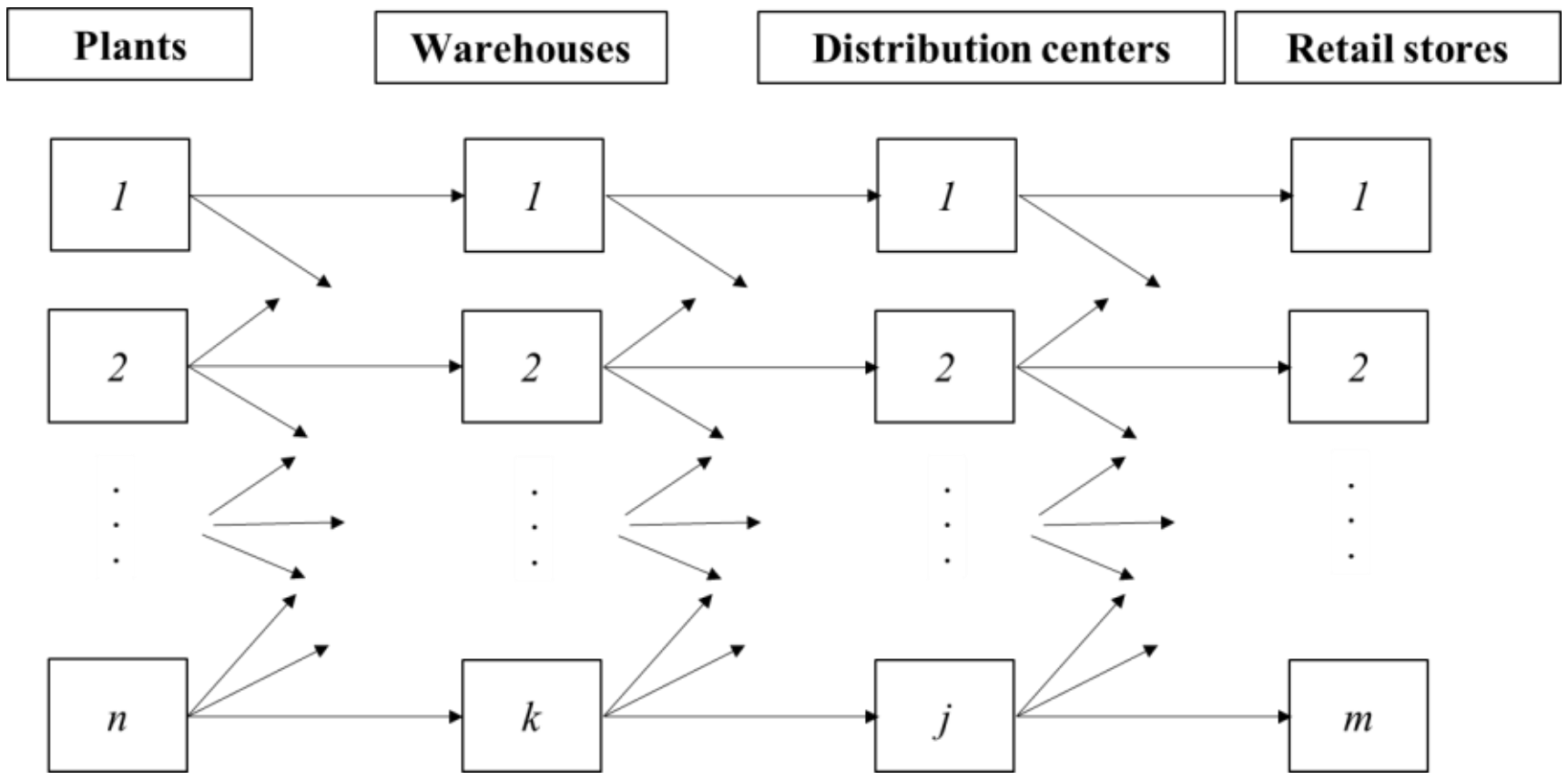

Figure1. A general topology of the 4L-FLP supply chain

Each arc (*a, b*) in the arc-set *A* is associated with a transportation cost, $c_{ab}$, which represents the expense of shipping a single bundle of products along that arc. Additionally, serving a retail store *m* is expected to generate a benefit, $R_m$. The 4L-FLP aims to optimize the selection of markets, plants, warehouses, and distribution centers to maximize the total profit of serving these selected markets. However, this selection process is subject to several constraints. Each selected facility incurs a one-time fixed cost, including the costs of operating a market, plant, warehouse, and distribution center. Resource limitations make it necessary to limit upper bounds on the number of selected stores, plants, warehouses, and distribution centers. These fixed costs can be viewed as annual operation expenses, while the single bundle of products represents a typical shipment to a store (Vanovermeire & Sörensen 2014; Ramshani. 2019; Myung & Yu 2020; Cao 2022; Tang 2023). The 4L-FLP is a crucial focus in strategic decisions for supply chain design, as highlighted in various studies (Lee & Whang 1999; Khalifehzadeh 2015; Saldanha-da-Gama 2022).

Research on *k*L-FLP often fails to address real-world complexities. Significant knowledge gaps include: 1) Assuming all customers must be served and plants are pre-identified, rather than considering the selection of markets, plants, warehouses, and distribution centers. 2) Neglecting the selection of facility numbers and upper bounds due to resource constraints. 3) Focusing on one-time fixed costs for warehouses and distribution centers, ignoring fixed costs for selecting all facilities. 4) Using mathematical formulas with numerous variables, making them impractical for large-scale problems and omitting important industrial issues. 5) Lacking sophisticated meta-heuristics to solve large-scale real kL-FLP, effectively. A detailed discussion of these knowledge gaps is given in section 2.5 below.

To address the limitations in the literature, our study employs a carefully crafted experimental framework to ensure rigorous validation of the proposed algorithms for the 4L-FLP. We establish precise objectives of 4L-FLP and develop specific, verifiable hypotheses to guide our research. We generate a diverse array of 45 problem instances to comprehensively evaluate algorithms' performance and carefully design to encompass a wide spectrum of sizes, structures, and complexities. Our experimental protocol follows strict standards, incorporating multiple independent runs with varying random seeds and employing systematic hyperparameter optimization to guarantee equitable algorithm comparison. We utilize key performance indicators, including Best-Found Solution (BFS) and Time to Best (TB), to accurately assess the solution and its efficiency. To prove the reliability of our findings, we conducted extensive sensitivity analyses, scrutinizing potential biases and limitations within our experimental design through statistical methods, including significance testing, non-parametric tests, and effect size calculations. In the interest of scientific transparency and reproducibility, we have made our complete dataset, encompassing all 50 instances and their corresponding solutions, publicly accessible.

To address limitations in the literature, this paper addresses the following issues:

(1) We tackle a 4L-FLP in the context of managing the daily supply chain where the selection of markets (retail stores or market zones), plants, warehouses, and distribution centers are made from available options. This approach takes a more holistic and realistic view in contrast to the existing literature, which typically assumes all customers are fixed and must be served, and all plants are known and fixed.

(2) To reflect real-world constraints, we impose upper bounds on selected facilities at each level.

(3) We incorporate *one-time* fixed costs associated with selecting a facility at each level of the supply chain.

(4) We propose a novel integer programming formula based on the quadratic assignment problem. This model is also arc-based and significantly reduces the number of variables compared to linear integer programs. Additionally, we provide the necessary mathematical results to enable the fast implementation of heuristics.

(5) We develop and experimentally test two solution procedures: a multi-start greedy and a multi-start tabu.

The rest of this paper is organized as follows. We begin by presenting a comprehensive literature review and summarizing knowledge gaps, providing a foundation for our research. Next, we introduce a novel mathematical formula for the 4L-FLP from a variation of the presented problem, followed by a section that derives several critical results necessary for developing heuristics. We propose two heuristics: a multi-start greedy local, and a multi-start tabu search. The performance of these heuristics is evaluated through a series of experiments presented in the subsequent section. We conduct comprehensive statistical analyses to assess the reliability of proposed heuristics based on the computational results. We conclude the paper by discussing the managerial implications, summarizing our findings, and outlining potential directions for anticipated research.

## 2. Literature

Single-level 4L-FLPs have been studied in various contexts. Their model topologies often resemble the graph in Figure 1, featuring *suppliers*, *plants*, *distribution centers*, and *customers*, and often using different names, depending on the application. For example, see (Majumdar 2023; Diglio 2024; Lin 2024; Tapia-Úbeda 2024). To avoid confusion, we adopt consistent terminology throughout this paper, categorizing facilities as *plants*, *warehouses*, *distribution centers*, and *markets* (customer zones or retail stores), as depicted in Figure 1. Contrary to existing literature, all customers are typically assumed to be fixed and must be served. All plants are known and fixed. See examples, (Aardal 1999; Melo 2006; Melo 2009; Latha Shankar 2013; Melo 2014; Ozgen & Gulsun 2014; Khalifehzadeh 2015; Byrka 2016; Eskandarpour 2017; Fattahi 2017; Zokaee 2017; Cortinhal 2019; Shoja. 2019; Kumar 2020; Belieres 2021; Vishnu 2021; Guo 2022; Majumdar 2023). Our approach distinguishes itself by considering a more comprehensive and realistic scenario. For comprehensive reviews of hierarchical facility location problems, refer to survey papers (Ortiz-Astorquiza 2018; Farahani 2019; Kumar 2020; Dukkanci 2024).

Most research on facility location problems concentrated on the two-level case, as discussed in studies (Kratica. 2014; Gendron 2017; Malik 2022; Gendron 2023; Sluijk 2023). However, the more general *k*-level facility location problem (*k*L-FLP) has also received attention from the researchers. Over the past two decades, several comprehensive surveys have summarized state-of-the-art in this field, including works by (Klose & Drexl 2005; Şahin & Süral 2007; Melo 2009; Farahani 2014; Ortiz-Astorquiza 2018; Farahani 2019; Kumar 2020; Dukkanci 2024). These surveys reveal insights into the development of *k*L-FLP research over the last three decades.

An extant literature review of *k*L-FLP suggests several shortcomings in multi-level facility location research. For a comprehensive review of hierarchical facility location problems, we refer to (Ortiz-Astorquiza 2018; Farahani 2019; Kumar 2020; Dukkanci 2024). For example, Marianov and Eisel (2024) considered a selected survey of location analysis in the last 50 years, including hierarchical facility location. However, due to the lack of published research, they focused on 2 and 3-echelon models, and due to the complexity of the problem, they emphasized the need for developing powerful techniques that profoundly deal with large-scale problems. Drezner & Eiselt (2024) also surveyed competitive location models. This model tries to find the best locations for facilities among existing and competing sites. The authors also emphasized the need for heuristics, especially meta-heuristics, for realistic multiple-facility location models.

Supply chain structure involves a series of interconnected decisions: selecting the number of retail stores, plants, warehouses, and distribution centers (Farahani 2019; Kumar 2020; Guo 2022; Kidd 2024; Lin 2024). Despite the complexity of this issue, existing literature suggests that researchers often need to pay more attention to incorporating all important features into their models (Farahani 2019; Kidd 2024). For example, Farahani (2019), in a survey of OR models on facility locations, provided a series of important features that needed to be included in these models (explained below). Kidd (2024) concentrated on a two-echelon supply chain network design facility location and proposed several realistic features in their model, including delivery date and flexible supply chain network design. Other notable exceptions that include some realistic features are considered in (Baumol & Wolfe 1958; Ortiz-Astorquiza 2017, 2018, 2019; Kidd 2024). Hajipour (2016) has gone a step further by incorporating the minimization of the number of facilities into the objective function of their model. This discussion highlights the need for a more comprehensive approach that considers the intricacies in relationships between these components of the supply chain structure.

Our literature review comprises mathematical formulas, including facility selection and fixed costs, upper bounds for the number of facilities, heuristics, and knowledge gaps.

### *2.1.Formulas*

Mathematical formulas for the *k*L-FLP typically involve variations of integer programs, including integer and mixed-integer linear programs, and are solved using simple greedy or meta-heuristics. In a recent survey, Ortiz-Astorquiza (2018) provided an overview of the formula, algorithms, and applications of *k*L-FLP, highlighting two prominent families of MIP formulas: arc-based and path-based models. While many authors have formulated *k*L-FLP and used heuristics to solve them, most formulas are variations of these two models. Here, we highlight some of the original papers in this area. For instance, Ortiz-Astorquiza (2018) presented an arc-based formulation that includes many factors, which are also considered in this paper. In contrast, our model in this study is an arc-based quadratic model. It can be considered closest to their model in terms of included factors. Other notable contributions include Melo (2006), who gave an MIP formula of the hierarchical facility location problem in the context of relocation. Şahin and Süral (2007) presented a linear assignment-based formula of the hierarchical facility location, a simplified version of the arc-based model, and applied simulated annealing to small, generated problems. Gendron and Semet (2009) discussed path-based and arc-based formulas of hierarchical facility location. They showed that linear programming relaxation of the path-based model provides a better bound than the LP relaxation of the arc-based model. Gendron (2017) recently compared six MILP formulas with a single assignment constraint for two-level facilities. Karatas and Eriskin (2023) considered hierarchical location and sizing problems and presented linear and piece-wise linear models. Although the above approaches are linear, when applied to realistic large-scale problems, the number of variables is exceedingly high, making it impractical to solve the problems. Our formula here is also arc-based; however, it is quadratic in the objective function based on the quadratic assignment problem (QAP). This formula significantly reduces the number of variables compared to the linear formulas. Many heuristics and meta-heuristics developed in the last three decades for QAP may be appropriately tuned for immediate use.

### *2.2.Facility Selection and Fixed Costs*

A common assumption in the literature is that markets and plants are pre-identified and must be served. In contrast, the locations of warehouses and distribution centers must be determined, and one-time fixed costs are incurred for opening them. This situation is considered in the following published papers: Kaufman (1977), Gendron and Semet (2009), Latha Shankar (2013), Ekici (2013), Ozgen and Gulsun (2014), Ortiz-Astorquiza (2017, 2018, 2019), Eskandarpour (2017), da Silveira Farias (2017), Fattahi (2017), Ou-Yang and Ansari (2017), Zokaee (2017), Cortinhal (2019), Lai (2019), Shoja (2019), Ruvalcaba-Sandoval (2021), Vishnu (2021), Kumar and Kumar. (2024), Menezes (2024), Vishnu (2021). The fixed costs can be interpreted as yearly operating expenses, and single products can represent cases where a bundle of products is shipped to a store (Vanovermeire & Sörensen 2014; Myung & Yu 2020; Cao 2022; Ridderstedt & Nilsson 2023; Tang 2023). However, in many supply chain designs, it is essential to select markets as well as network facilities to serve those markets (Ortiz-Astorquiza 2017, 2018, 2019; Wang 2022). Therefore, this study focuses on optimizing the selection of a subset of potential markets, plants, warehouses, and distribution centers to maximize the total profit generated from serving the chosen markets.

### *2.3.Number of Facilities*

In a recent survey on the integration of planning in supply chains, Kumar (2020), and Marianov and Eiselt (2024), highlight the importance of determining the facility numbers at each level of the 4L-FLP. However, despite its relevance, research on facility selection has primarily overlooked the consideration of imposing an upper bound on the maximum number of facilities. In contrast, early work by Baumol and Wolfe (1958) formulated a mathematical model to identify the number and locations for warehouses. More recently, Latha Shankar (2013) developed a multi-objective optimization model for single-product 4-echelon supply chain architecture, which included identifying the required number of warehouses. In several papers, Ortiz-Astorquiza (2017, 2018, 2019) and Shavarani (2018) addressed the *k*L-FLP problem, restricting the number of facilities at each level to an upper bound. Menezes (2024) recently included a penalty term in the objective function that reduces the number of plant facilities. Furthermore, Hajipour (2016) proposed a multi-objective, and multi-layer facility location-allocation model. This aims to determine the optimal number of facilities and service allocation at each layer, with the number of facilities to be chosen included as a goal in the objective. In this paper, we limit the selection of the number of markets, plants, warehouses, and distribution centers to upper bounds.

### *2.4.Heuristics*

Due to the complex nature of the *k*L-FLP, many researchers have suggested the need for developing heuristics that can be applied to the biggest problems. Earlier, Şahin and Süral (2007) employed simulated annealing to tackle hierarchical facility location problems, achieving good average solutions. However, their experiment is based on very small problems. Notably, their largest problem instance involved 150 sites. Their study stands out for its pioneering application of meta-heuristic approaches to hierarchical problems. Surprisingly, we have not seen the use of sophisticated neighborhood search methods, or meta-heuristics, to solve realistic large-scale hierarchical location problems. In a recent survey of OR models in service facility locations, Farahani (2019) showcased the importance of meta-heuristic techniques in optimization modeling, particularly when dealing with complex problem modeling or large-scale optimization challenges. The authors emphasized that, despite their potential, meta-heuristics have been underutilized compared to heuristics, creating a gap in Operations Research. Latha Shankar (2013) also addressed multi-objective optimization for a single-product, four-echelon supply chain architecture. To optimize two objectives simultaneously, they used a four-echelon network model. They solved this problem using a swarm intelligence-based Multi-objective Hybrid Particle Swarm Optimization. Melo (2014) modeled the problem as a large-scale mixed-integer linear program and then developed a two-phase heuristic approach to obtain high-quality, feasible solutions. Khalifehzadeh (2015) introduced a heuristic algorithm called Swarm Optimization. Mortazavi (2015) integrated agent-based simulation techniques with reinforcement learning to model a four-echelon supply chain facing non-stationary customer demands. This approach enabled the authors to capture the supply chain dynamics, making informed decisions effectively. Hajipour (2016) proposed a Pareto-based multi-objective metaheuristic. The approach combines multi-objective vibration damping optimization and multi-objective harmony search algorithms to find and analyze Pareto's optimal solutions. Eskandarpour (2017) observed that a wide range of solution techniques, including exact and approximate approaches, have been employed to tackle facility location and supply chain network design problems. They note that the Large Neighborhood Search (LNS) technique, despite proven efficiency and flexibility in solving complex combined optimization problems, has been largely overlooked. Da Silveira Farias (2017) proposed a simple heuristic method for solving strategic supply chain design with four layers in a Brazilian tire company.

Ou-Yang and Ansari (2017) tackled the hierarchical facility location problem, as a complex optimization challenge that involves strategically positioning facilities to serve lower-tier facilities efficiently. They developed mixed-integer programming, which incorporates flow capacity constraints. They also proposed a hybrid approach combining Particle Swarm Optimization and Tabu Search as a solution procedure. The objective is to minimize the total demand-weighted distance traveled, facility operating costs, and flow assignment costs. Cortinhal (2019) presented a MILP to optimize multi-stage supply chain network design. The model determines plant and warehouse locations, supplier and transportation mode selection, and product distribution while minimizing costs and meeting customer service expectations. Janjevic (2019) developed a nonlinear optimization model for last-mile distribution, integrating collection-and-delivery point location decisions with demand pattern shifts. They introduced routing with cost approximation formulae and a heuristic solution method to enable scalability. The approach is validated through a real case study with a leading Brazilian e-commerce company, demonstrating its practical applicability and industry impact. Shoja (2019) probed the advantages of flexible supply chain networks, which have proven effective in manufacturing and service industries. The authors recognized that delivery modes can significantly impact overall supply chain efficiency. To address this, they proposed a mixed-integer linear programming (MILP) model for designing a multi-product, four-stage flexible supply chain network in a solid transportation environment. The model is solved using ten top meta-heuristic algorithms. Wu (2019) tackled the challenge of solving complex supply chain optimization problems, which often exceed the capabilities of commercial-led solvers. To overcome this limitation, a novel regression and extrapolation-guided optimization method capable of efficiently solving medium- and large-sized problem instances was proposed. Belieres (2021) focused on the operations of a third-party logistics (3PL) service provider in supply chain management for a French restaurant chain. The authors developed a network reduction matheuristic, inspired by the Dynamic Discretization Discovery algorithm, optimizing 3PL provider's operations and improving overall supply chain efficiency. Kumar and Kumar (2024) proposed mixed-integer linear programming (MILP) for designing an uncertain supply chain network that minimizes overall costs while considering carbon emissions and social factors. A greedy-based heuristic for solving larger instances is proposed, and sensitivity analysis is conducted to explore the impact of various parameters. Ortiz-Astorquiza (2019) introduced an exact algorithm based on a Benders decomposition to solve multi-level uncapacitated p-location problems. They efficiently generate Pareto-optimal cuts by leveraging the network flow structure. Computational experiments on benchmark instances with up to 3,000 customers, 250 facilities, and four levels demonstrate the algorithm's efficiency and effectiveness.

Recent studies have focused on various aspects of supply chain network design. For instance, Ruvalcaba-Sandoval et al. (2021) proposed a four-level supply chain network design model that determines the number, locations, and capacities of factories and warehouses, as well as transportation between different sites. They developed a MILP model and two matheuristic algorithms to solve the problem. Similarly, Belieres (2021) considered a 4-echelon network design for restaurant supply chains, proposing a MILP and a matheuristic algorithm as a solution-finder. Janjevic (2019) integrated collection-and-delivery points (CDPs) into the design of multi-echelon distribution networks, developing a non-linear optimization model and a heuristic solution method for Brazilian e-commerce last-mile delivery. Additionally, several approximation algorithms with guaranteed bounds have been proposed in the literature, including works by (Aardal 1996) Aardal (1999) Zhang and Ye (2002); Kantor and Peleg (2009) Drexl (2011) Li (2013) Ortiz-Astorquiza (2017).

Overall, the literature emphasizes the need to develop meta-heuristics for multi-echelon facility locations. For recent examples, see the following papers; all emphasize such need (Aardal (1999)). To address this gap, this paper provides two heuristics (1) a multi-start greedy heuristic and (2) a multi-start tabu search.

The greedy heuristic is based on the single flip of a facility, using Proposition 1 (Appendix A). Generally, starting solutions in multi-start searches randomly select a starting solution, which typically leads to poor solutions, Drezner (2024). However, we adopt a mechanism usually used in sequencing problems, such as a traveling salesman creating diverse solutions, known as the *r*-Opt strategy. We also use the *r*-Opt strategy as the search is in progress. This creates an opportunity to explore diverse areas of solution space. The approach has been applied to several integer programs with considerable success; see, for example (Alidaee & Wang 2017; Wang & Alidaee 2019; Wang & Alidaee 2023). The greedy heuristic is further used as an improvement stage in the tabu search.

### *2.5.Knowledge Gaps and Research Motivation*

In a recent comprehensive survey of OR models in facility locations, Farahani (2019) highlights the typical input data required for *k*L-FLP, including the locations of current and potential facilities, facility costs, capacity, distance, response time, and service level. This review also reveals significant gaps between research and real-world scenarios that beg serious attention. The authors specifically address the following gaps in the literature, and we include these issues in our model.

1. *Market, plant, warehouse, and distribution center selection:* Researchers often assume that all customers must be served and all plants are pre-identified, but in strategic supply chain management, it is crucial to consider site and total function selection for markets, plants, warehouses, and distribution centers.
2. *Facility number selection:* Researchers would do better to pay more attention to the selection of the number of facilities, despite the importance of considering upper bounds due to resource constraints.
3. *Fixed costs:* Most researchers only consider one-time fixed costs for selecting a warehouse and distribution center, whereas in many cases, such as entering new markets, it is necessary to consider the fixed costs of selecting all facilities.
4. *Formula limitations:* Existing mathematical formulas are often MIP models with numerous variables, making them inapplicable for large-scale real problems. These formulas also omit crucial industrial issues, highlighting the need for new models incorporating key industrial concerns while reducing variable numbers.
5. *Solution procedure limitations:* The literature needs more sophisticated meta-heuristics to solve large-scale real *k*L-FLP effectively. Therefore, it is essential to develop procedures to provide reasonable solutions for large-scale problems.

We explain some of the above shortcomings in (i)-(iv) below.

(i) The existing literature on facility locations mostly assumes that the number of retail stores (customers or market zones) is fixed, all must be served, and the locations of plants are also predetermined. This is evident in numerous studies, such as (Aardal *et al.* 1999; Melo *et al.* 2006; Melo *et al.* 2009, 2014; Ozgen & Gulsun 2014; Khalifehzadeh *et al.* 2015; Byrka *et al.* 2016; Eskandarpour *et al.* 2017; Zokaee *et al.* 2017; Cortinhal *et al.* 2019; Shoja *et al.* 2019; Kumar *et al.* 2020; Belieres *et al.* 2021; Vishnu *et al.* 2021; Guo *et al.* 2022; Yan *et al.* 2022; Majumdar *et al.* 2023). However, in real situations, it is top priority to determine which markets to serve, where to locate plants, and where to establish warehouses and distribution centers.

(ii) In all these studies, one-time fixed costs are only considered for warehouses and/or distribution centers, as seen in (Kaufman 1977; Gendron & Semet 2009; Ekici 2013; Ozgen & Gulsun 2014; Eskandarpour 2017; Fattahi 2017; Ortiz-Astorquiza 2017; Zokaee 2017; Ortiz-Astorquiza 2018; Cortinhal 2019; Lai 2019; Ortiz-Astorquiza 2019; Shoja 2019; Vishnu 2021). However, as already explained, it is necessary to determine which markets to serve, where to locate plants, and where to establish warehouses and distribution centers in real-world scenarios. In these situations, selecting a facility ensures a one-time fixed cost. A recent study by (Sebatjane & Adetunji 2024) also highlights the importance of considering fixed costs incurred when processing batches of products sent to retailers.

(iii) The solution methodologies for the *k*L-FLP typically involve variations of mathematical formulations, such as integer and mixed-integer linear programs, which are then solved using simple heuristics or meta-heuristics. These formulas can be categorized into *assignment-based*, *arc-based*, and *path-based* models, as seen in the works of (Melo 2006; Şahin & Süral 2007; Gendron & Semet 2009; Kratica 2014; Ortiz-Astorquiza 2017, 2018; Farahani 2019; Ortiz-Astorquiza 2019; Karatas & Eriskin 2023). Several authors have studied the 4-level facility location problem (4L-FLP) in various contexts, (Mohammadi *et al.* 2017; Cortinhal 2019; Wu 2019; Belieres 2021; Guo. 2022; Badejo & Ierapetritou 2023; Majumdar 2023; Sebatjane & Adetunji 2024). In a recent survey of reinforcement learning on logistics and supply chain management, Yan (2022) discuss the 4-level supply chain location analysis. Similar to the general *k*L-FLP, formulations for the 4L-FLP can be either arc- or path-based. Saldanha-da-Gama (2022) emphasized the need for studies on more complex and challenging problems, leading to development of comprehensive mathematical models, particularly for large-scale problems. Here, we propose a novel arc-based mathematical programming formulation based on a variation of quadratic assignment problems (QAP), significantly reducing the number of variables. Furthermore, the vast literature on large-scale applications of heuristics and meta-heuristics for QAP can be leveraged to solve 4L-FLPs.

(iv) Most studies focus on developing algorithms for two-level facility locations, see for recent works Saldanha-da-Gama (2022). However, several simple and meta-heuristics have also been applied to *k*L-FLPs, as seen in studies by (Şahin & Süral 2007; Latha Shankar 2013; Melo 2014; Khalifehzadeh 2015; Mortazavi 2015; Hajipour 2016; Eskandarpour 2017; Fattahi 2017; Ou-Yang & Ansari 2017; Cortinhal 2019; Janjevic 2019; Shoja 2019; Wu 2019; Belieres 2021; Guo 2022; Kumar & Kumar 2024). These heuristics mostly concentrate on small to medium size problems. Furthermore, they consider only some of the realistic features in their models. Interestingly, all of these publications emphasize the need to develop heuristics for realistic large-scale problems. In a comprehensive survey of hierarchical facility location problems, Farahani (2019) highlighted that researchers often resort to heuristic or meta-heuristic techniques due to the complexity of modeling or the size of the problem. They note that most heuristics used for multi-echelon facility location problems are simple heuristics. They suggested that more attention should be focused on the need for the development of meta-heuristic techniques. This gap needs to be addressed from an operations research (OR) perspective. This sentiment is echoed by Cortinhal (2019), who emphasize the need for developing heuristic procedures for multi-echelon facility locations. This is also evident from a recent article by (Kumar & Kumar 2024), where they use a simple greedy solution procedure. Da Silveira Farias (2017) pointed out that solution procedures are mainly based on heuristic techniques, which have severe limitations for large-scale problems. The large-scale nature of these problems also makes commercial solvers' branch-and-bound methods unfavorable, as noted by Wu (2019), who proposes a regression and expanded-thinking guided method as an alternative solution approach. Eskandarpour (2021) showed a local search procedure as it were a supply chain network design problem. Guo (2022) also emphasized the need for further development of large-scale solution methods. Ortiz-Astorquiza (2017) and Ortiz-Astorquiza (2019) have recently proposed approximation and exact algorithms for $k$L-FLPs. Several approximation algorithms with guaranteed bounds have also been proposed; refer to (Ortiz-Astorquiza 2017, 2018; Farahani 2019) for a detailed analysis of these results.

Ortiz-Astorquiza (2018) noted that recent variants of the $k$L-FLP have emerged, allowing planners to choose between incurring a penalty and serving all customers. This approach results in customer selection, which shares similarities with our customer selection considerations in certain aspects of this study.

## 3. A QAP Formula and Solution Procedure

In this section, we present a mathematical formula based on a quadratic assignment problem (QAP) for 4L-FLP. Then, a greedy multi-start local search and a sophisticated multi-start tabu search.

The 4L-FLP can be mathematically formulated as a complex bi-quadratic programming problem, expressed as follows.

$S$ Set of potential retail stores (customer zones or markets) to be served, indexed by $m$

$P$ Potential plants, indexed by $n$

$W$ Potential warehouse locations, indexed by $k$

$D$ Potential distribution center locations, indexed by $j$

$P_m$ Plants eligible to deliver products to a retail store $m$

$D_m$ Set of distribution centers eligible to deliver products to a retail store $m$

$G= (V, A)$ is a graph with nodes $V$ and directed arcs $A$, where:

$V = P \cup W \cup D \cup S$

$A = PW_e \cup WD_e \cup DS_e$

$PW_e = \{e = (n, k): n \in P, k \in W\}$, All arcs between plants and warehouses

$WD_e = \{e = (k, j): k \in W, j \in D\}$, All arcs between warehouses and distribution centers

$DS_e = \{e = (j, m): j \in D, m \in S\}$, All arcs between distribution centers and retail stores

$c_{ab}$ Costs of transporting a bundle of products along an arc $(a, b) \in A$

$fp_n$ Costs of opening a plant *n* in *P*

$fw_k$ Costs of opening a warehouse *k* in *W*

$fd_j$ Costs of opening a distribution center *j* in *D*

$fs_m$ Costs of opening a retail store *m* in *S*

*up* Upper bound for new plants to be opened

*uw* Upper bound for new warehouses

*ud* Upper bound for new distribution centers

*us* Upper bound for new retail stores

$R_m$ Total revenue if retail store (market) *m* is served

Decision Variables:

$z_{pn}$ is 1 if a new plant $n \in P$ is opened, and 0 otherwise

$z_{wk}$ is 1 if the new warehouse $k \in W$ is opened, and 0 otherwise

$z_{dj}$ is 1 if new distribution center $j \in D$ is opened, and 0 otherwise

$z_{sm}$ is 1 if new retail store (market) $m \in S$ is served, and 0 otherwise

$x_{mnk}$ is 1 if an operational plant *n* ships a bundle of products to an open retail store *m* via an open warehouse *k*, and 0 otherwise.

$y_{jm}$ is 1 if an opened distribution center *j* delivers a bundle of products to an opened retail store *m*, and 0 otherwise

*TR* Total revenue

*TC* Total fixed and transport costs

*TP* Total net profit

$$TR = \sum_{m \in S} R_m z_{sm}$$

$$TC = \sum_{m \in S} (fs_m) z_{sm} + \sum_{n \in P} (fp_m) z_{pn} + \sum_{k \in W} (fw_k) z_{wk} + \sum_{j \in D} (fd_j) z_{dj} +$$

$$\sum_{m \in S} \sum_{n \in P_m} \sum_{(n,k) \in PW_e} c_{nk} x_{mnk} + \sum_{m \in S} \sum_{n \in P_m} \sum_{(n,k) \in PW_e} \sum_{(k,j) \in WD_e} \sum_{j \in D_m} c_{kj} x_{mnk} y_{jm} + \sum_{(j,m) \in DS_e} c_{jm} y_{jm}$$

*Max* $TP = TR - TC$ (1)

s.t. $\sum_{j \in D_m} y_{jm} = z_{sm}, \ \forall m \in S$ (2)

$$\sum_{n\in P_m}\sum_{(n,k)\in PW_e} x_{mnk} = z_{sm},\ \forall m \in S \quad (3)$$

$$y_{jm} \le z_{dj},\ \forall m \in S, j \in D_m \quad (4)$$

$$x_{mnk} \le z_{pn},\ \forall m \in S, n \in P_m, (n,k) \in PW_e \quad (5)$$

$$x_{mnk} \le z_{wk},\ \forall m \in S, n \in P_m, (n,k) \in PW_e \quad (6)$$

$$\sum_{n\in P} z_{pn} \le up \quad (7)$$

$$\sum_{k\in W} z_{wk} \le uw \quad (8)$$

$$\sum_{j\in D} z_{dj} \le ud \quad (9)$$

$$\sum_{m\in S} z_{sm} \le us \quad (10)$$

$$x_{mnk}, y_{jm}, z_{sm}, z_{pn}, z_{wk}, z_{dj}\ \ Binary,\ \forall\ m \in S, n \in P, k \in W, j \in D \quad (11)$$

The objective function (1) maximizes profit from strategically placing plants, warehouses, and distribution centers to efficiently serve all opened retail stores. Constraint (2) ensures that if a store is opened, it must be serviced through a designated distribution center. Similarly, constraint (3) ensures an eligible plant provides products to any opened store. Constraints (4)-(6) guarantee that if a store is served, distribution centers, plants, and warehouses are opened to facilitate the service. Constraints (7)-(10) uphold the upper bounds for plants, warehouses, distribution centers, and retail stores. Constraint (11) ensures the satisfaction of binary variables. It is noteworthy that the objective function involves a biquadratic term. Specifically, when a store (*m*) is opened and served through a distribution center (*j*), the associated arc (*j*, *m*) is activated $y_{jm} = 1$. Similarly, for a store (*m*) to be served, a corresponding plant (*n*) and warehouse (*k*) must be opened $x_{mnk} = 1$. Consequently, the biquadratic term $y_{jm}x_{mnk} = 1$ in the objective function ensures the activation of the arc (*k*, *j*) warehouse (*k*) and distribution center (*j*) as well.

Solving 1-11, we present heuristics that can be applied to large-scale problems. To implement the heuristics, we need several results. Proposition 1 in Appendix A serves this purpose. We give our heuristics in Appendix B.

Note that the objective function has a biquadratic term. This captures the complex relationships between the facilities in the supply chain. Specifically, when a retail store (*m*) is opened and served through a distribution center (*j*), the corresponding arc (*j*, *m*) is activated, indicated by $y_{jm} = 1$. Furthermore, for a store (*m*) to be served, a corresponding plant (*n*) and warehouse (*k*) must be operational, which is denoted by $x_{mnk} = 1$. The biquadratic term $y_{jm}x_{mnk} = 1$ in the objective function ensures that the arc (*k*, *j*) between the warehouse (*k*) and the distribution center (*j*) is also activated, thereby maintaining the integrity of the supply chain.

A feasible schedule for a retail store *m* corresponds to a feasible path (*n*, *k*, *j*, *m*) in the supply chain graph *G*, which defines the sequence of facilities that serve store *m*. The topology of the 4-level supply chain, in Figure 1, exhibits the relationships between different levels of facilities. We propose a multi-start greedy local search, **Heuristic 1**, which incorporates an embedded *r*-Opt sequence diversification strategy, as detailed in Appendix B. For a comprehensive discussion of r-Opt strategies within heuristics, refer to (Alidaee & Wang, 2017). Heuristic 1 is used as an improvement strategy in the multi-start tabu search, **Heuristic 2,** as detailed in Appendix B.

## 4. Experimental Design and Results

### *4.1 Data and Experimental Design*

There is no benchmark available for the problems we consider in this paper. The only benchmark that has some characteristics of our problems is given by Ortiz-Astorquiza (2019). However, among them are only two with four-level facilities. These two problems also miss some data that cannot be tuned for use in our computational experiment. Thus, we randomly generated problem instances and solved them using local search, Heuristic 1, and multi-start tabu search, a key component of Heuristic 2. All algorithms were implemented in Fortran and executed in order on the core of an Intel Xeon Quad-core E5420 Harpertown processor equipped with a 2.5 GHz CPU with 8 GB of memory. The parameters used to generate the problem instances are detailed in Table A1; see Appendix D.

Several key considerations emerged during the development of our heuristics. One crucial aspect was the incorporation of sequences in the heuristics, a strategy that has proven successful in other combined optimization problems (Alidaee & Wang, 2017; Wang & Alidaee, 2019, 2023). We integrated the *r-Opt* strategy into Heuristic 1, initially experimenting with 1-Opt, 2-Opt, 3-Opt, and a limited version of 4-Opt. Our results showed that 2-Opt, 3-Opt, and limited 4-Opt were the most effective approaches, with limited 4-Opt ultimately chosen for all heuristics. Notably, the required time to reach the top solution decreased by about 67%, on average. The embedded sequences within the heuristics created the chance to explore a larger area of the solution space, generating multiple starting points as the search progressed, each yielding distinct solutions depending on the sequence used (see Figure 3).

The placement of the tabu strategy in a heuristic significantly impacts its performance. Our experiments revealed that applying the tabu strategy to variables at different levels of the supply chain yields varying results. Implementing it in the selection of retail stores led to improvements, but applying it to the selection of plants, warehouses, and distribution centers often obscured the outcomes. Consequently, optimizing its benefits, we strategically integrated the tabu strategy only in the retail store selection component of Heuristics 2 and 3.

The primary objective of meta-heuristics like tabu search is to escape local optimality and explore a broader solution space. When a local optimum is reached, a common strategy is to perturb the solution, a process known as "shaking," then re-optimize, or "baking." This explores new areas of the solution space. We tested this approach by employing different shaking and baking intensities, randomly choosing a number between L1 and L2, perturbing the solution, and solving problems of varying sizes (refer to Appendix D, Table A1 for the definition of parameters). In the experiments, the density of all matrices was 20%. These include matrices of m-p (a 0-1 matrix of m by n, 1 means retailer m can use plant 1, 0 otherwise), p-w (an n by k matrix of cost transporting a bundle of products from plant n to warehouse k), w-d (a k by j matrix of cost transporting a bundle of products from its warehouse k to the distribution center j), and d-m (a j by m matrix of costs transporting a bundle of products from distribution center j to retail store m). Similar density was used to choose several stores, plants, warehouses, and distribution centers. Similarly, the upper bounds for the several retailers, plants, warehouses, and distribution centers used to be 20%.

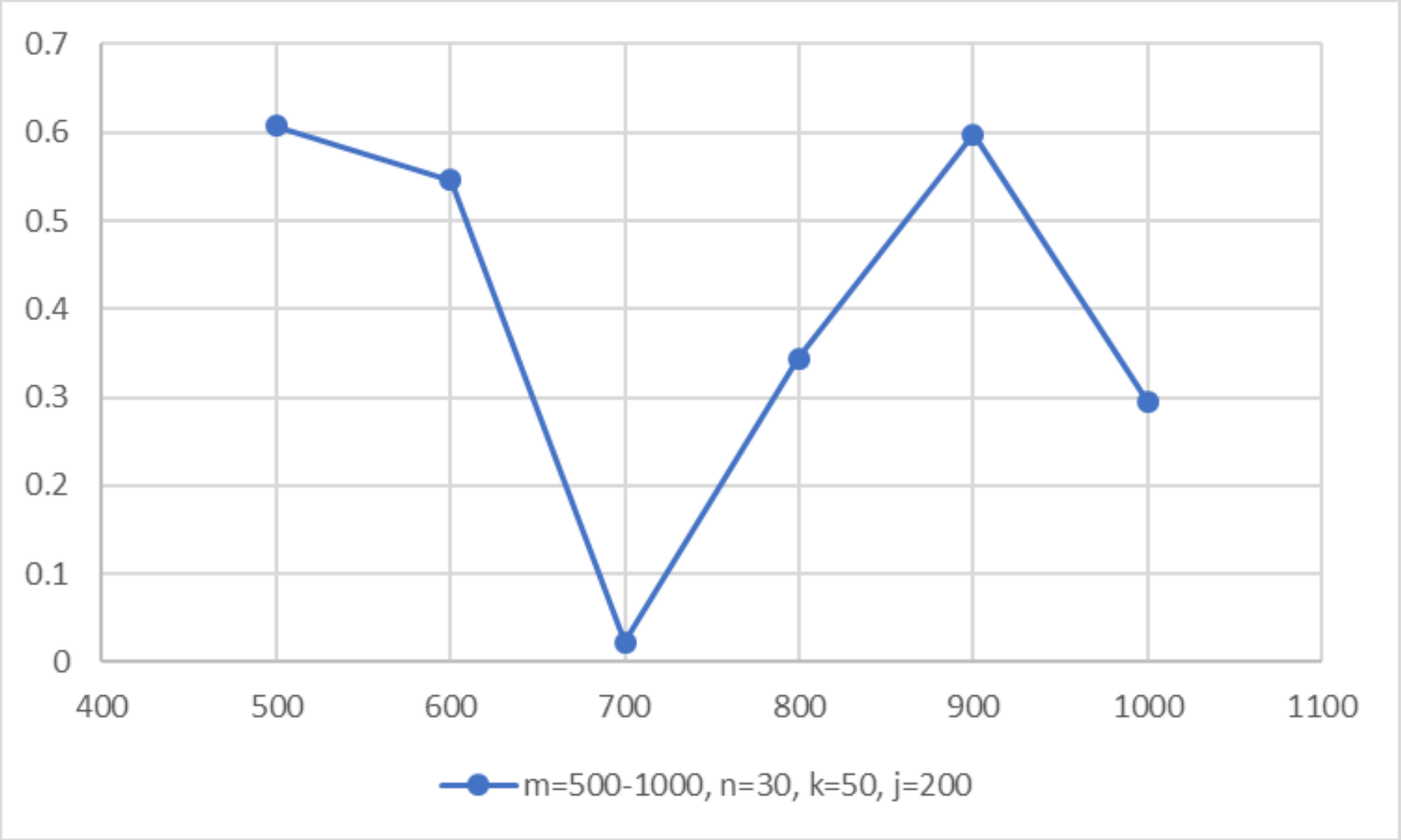

Figure 3. Percentage improvement in LS using a limited 4-Opt strategy. Time to best also improved by 67%.

Tabu tenure size is critical in the tabu search process, and its optimal value can significantly impact the algorithm's performance. We conducted experiments to determine the ideal tabu tenure size for our problem instances, characterized by (*m, n, k, j*) = (1000, 30, 50, 200). Our results showed that setting the tabu tenure to 2.5% of *m* yielded a reasonable outcome, and we adopted this setting for subsequent experimentation.

We conducted a comprehensive experimental study to test the performance of the local search and tabu computer-search algorithms on problems of varying sizes. Specifically, we generated instances with *m* (number of retailers) from 2000 to 6000, *n* (plants) equal to 30, *k* (warehouses) equal to 50, and *j* (distribution centers) equal to 150. We first used Heuristic 1 with and without *r-Opt* sequencing diversification. 10 instances were generated randomly and solved 15 times (runs), each time with different starting solutions. Table 1 represents the *BFS* and *TB* for each of the 10 generated instances. As can be seen from the tables, the *BFS* and *TB* for local search with sequence diversification significantly outperform when no sequence diversification is implemented. Then, we used the local search and tabu search (Heuristic 3), both with sequence diversification implemented, and solved 10 cases. The local search was used to solve the same instances 15 times (runs) with different starting solutions. The tabu search was used to solve the 10 instances given 2 seconds. However, within 2 seconds, we let the search continue until it stops. Table 2 represents the *BFS* and *TB* for each of the 10 generated instances by each algorithm. As seen in this table, the tabu search significantly outperforms the local search. Interestingly, the *TB* to reach the *BFS* is significantly small in comparison with local search.

[insert Tables 1 and 2 here]

### *4.2 Sensitivity Analyses of algorithms*

To test the reliability and sensitivity of the proposed algorithms, this study employs a comprehensive suite of statistical methods. The foundation of this approach involves conducting multiple independent runs with different random seeds to account for the random nature of the used optimization algorithms. Descriptive statistics, including the preset and expected light deviation of BFS improvements and TB reduction, provide an overview of the performance and consistency of the algorithms (LS and TS) across these runs. Table 3 shows that Local Search with sequence implementation improves BFS on all instances across 15 runs, with consistent improvements relative to the BFS value of Local Search without sequence, ranging between 0.68% and 0.81%, as the problem size increases. Meanwhile, Local Search with sequence significantly reduces TB on all instances across 15 runs, with a time reduction between 85% and 90%. The results of 15 runs show that Tabu Search with sequence outperformed Local Search with sequence across all instances. The consistent improvements, measured by *BFS_diff* values, range from 1.1% to 2.0%, while *TB_diff* values range from 73% to 80%. These findings indicate that the performance of Tabu Search with sequence is reliable across different instance sizes and the variation within instances of the same size.

[insert Table 3 here]

Hypothesis testing, such as t-Test or ANOVA, can be used to compare the different algorithms or parameter settings, determining if observed differences are statistically significant. In this study, ANOVA for one single factor and the t-Test (paired two samples for means) are used as statistical methods to compare group means. The t-Test calculates a t-value based on the difference between the two means relative to the variability whereas ANOVA uses an F-value derived from the difference of between-group variance to within-group variance (Maxwell *et al.* 2017). When the sample size is small (e.g., n = 10 in this study for the set of instances with the same size), the F-value is typically less significant than the t-value because ANOVA distributes total variance across multiple groups, increasing the degrees of freedom and reducing statistical power (Gelman 2005). This makes it harder to detect significant differences unless the effect size is large. The t-Test is more sensitive for two-group comparisons, while ANOVA is more appropriate for three or more groups, as it controls inflated error rates from multiple comparisons (Maxwell *et al.* 2017). Table 4 shows that there are significant differences between BFS value of Local Search with sequence and of Local Search without sequencing based on t-value. F-values of three instances sets (3000-30-50-150, 5000-30-50-150 and 6000-30-50-150) on BFS are statistically significant with p-value<=0.05 between two Local Search implementations. However, both F-value and t-value on TB are significant across all instances. These results are consistent with Table 3. Meanwhile, the F-values and t-values of BFS and TB show significant differences between Tabu Search with sequence and Local Search with sequence across all instances. These findings indicate that the performance of Tabu Search with sequence is statistically reliable across different instance sizes and the variation within instances of the same size.

[insert Table 4 here]

When data does not meet normal assumptions, which is common in optimization results, non-parametric tests like the Wilcoxon signed-rank or Mann-Whitney U are more effective. These tests evaluate whether there is a visible difference between the medians of two groups by ranking the absolute differences between paired observations while accounting for their signs. A p-value of 0.02 indicates a 2% probability that the observed difference occurred by random chance, providing strong evidence to reject the lesser or invalid hypothesis at the 5% significance level ($\alpha$ = 0.05). This suggests a statistically significant difference between the two paired samples, meaning the observed effect is unlikely due to randomness. As shown in Table 5, the BFS and TB values of Local Search with sequence differ significantly from those of Local Search without sequence. Similarly, the BFS and TB values of Tabu Search are statistically different from Local Search with sequence, demonstrating continuous improvement across algorithms. For instance, a p-value of 0.02 supports the conclusion that Tabu Search led to significant improvements over Local Search.

[insert Table 5 here]

## 5. Managerial Implications

As previously mentioned, the 4L-FLP model with the characteristics considered in this study is particularly well-suited for supply chain design, especially when companies are entering or expanding existing markets. Notably, our research is the first to consider all the aspects we have explained earlier.

We employed the heuristic algorithm based on Theorem 1, which involves swapping each pair of facilities. We opted for these simple exchanges to reduce computational time. However, exploring more sophisticated local improvement processes in future research would be worthwhile. There are different methods to embed sequences within heuristics, as demonstrated by Alidaee & Wang (2017), Wang & Alidaee (2019, 2023), and it would be valuable to compare these approaches to determine which one performs best for these problems. Note that different sequences can lead to distinct outcomes each time a new sequence is chosen. This effectively employs a multi-start strategy that explores larger areas of the solution space. The choices of L1, and L2 in the heuristics are crucial, and we experimented with different combinations for various problems to determine the best values for these parameters (see Table 1). Nevertheless, these combinations can also impact computational time, and further research is needed in this area. Finally, while we focused on 4-layer supply chains in this study, it is essential to investigate more layers in future research to reflect the complexity of real-world supply chain designs.

In today's rapidly changing business landscape, optimizing facility locations is essential for companies to stay ahead of the competition. Multilevel facility location is a strategic approach. Optimally, selecting the location of facilities minimizes costs, maximizes efficiency, and improves customer satisfaction. Businesses can lower transportation costs, improve delivery, and increase response to changing markets by optimizing facility locations. However, managing multilevel facility locations poses significant challenges, including optimizing multiple facilities across different levels, balancing competing objectives, and adapting to changing markets, changing customer needs, and various external factors. The right selection of the location of retail stores, distribution centers, warehouses, and suppliers has helped major retailers like Wal-Mart and Home Depot experience significant growth over the last two decades. A large customer base, diverse supply base, wide product variety, multiple distribution channels, and internationalization characterize their supply chains. For example, Wal-Mart globally serves 255 million customers weekly, operates some 10,619 retail units in 19 countries, and offers over 140,000 SKUs with thousands of suppliers. Operating in a highly competitive environment, retail supply chains must meet elevated customer expectations for price expectations, delivery, and service, contributing to their complexity. To overcome these challenges, businesses must leverage advanced analytics, optimization techniques, and digital technologies to unlock the full potential of multilevel facility location optimization, driving growth, profitability, and sustainability in their operations.

The retail supply chain is vital to connecting customers with vendors. When retail takes center stage, it receives individual customer orders, fulfilled either from inventory held at retail locations or through direct shipments from manufacturers or wholesalers. This approach allows retailers to avoid holding inventory, a strategy often employed by internet-based catalog retailers. Large companies and corporations may control the retail and wholesale tiers, while third-party logistics providers handle direct shipments. Additionally, retailers often establish close relationships with manufacturers by ordering products with their brand names.

The retail supply chain consists of multiple levels, wholesale, and retail. The manufacturing and wholesale levels comprise many facilities with varying degrees of partnership persistence. For example, the very large Amazon.com uses leased distribution centers strategically located near customer zones and tax-friendly areas while exploring options for relocation and expansion. In contrast, apparel retailers like Zara and Benetton achieve agility by contracting with numerous small manufacturers, which can quickly adapt to changing customer demand. Companies like Zepter International are interested in conducting network analyses and finding ideal placements for their distribution centers to optimize their supply chains. By optimizing facility location, companies can reduce transportation costs, improve service levels, and decrease overhead costs, ultimately increasing their competitiveness in the market. (Soshko *et al.* 2007)

## 6. Conclusion

The 4L-FLP model is well-suited for supply chain design, particularly for companies entering or expanding markets. This study is the first to consider all aspects previously explained. The research employed a heuristic algorithm based on Theorem 1, involving facility swaps to reduce computational time. Embedding sequences within heuristics improved solution quality and reduced computational time, consistent with previous studies. The choice of parameters L1, and L2 in the heuristics is crucial, affecting results and computational time. Local Search can be time-consuming, and n, k, and j values for a given m significantly impact results. Balancing computational time and desired results is necessary, and further research is needed in this area. The study focused on 4-layer supply chains, but investigating more layers in future research is essential to reflect real-world complexity. Optimizing facility locations is crucial for businesses to minimize costs, maximize efficiency, and improve customer satisfaction. The retail supply chain connects customers with vendors, involving multiple levels and complex relationships between manufacturers, wholesalers, and retailers. Companies can increase competitiveness by optimizing facility locations to reduce costs and improve service.

**Data Availability Statement:**

**Data is available to download: https://figshare.com/s/7ff4e276e02542f3ded1**

## Reference

Aardal, K., Chudak, F.A., Shmoys, D.B., 1999. A 3-approximation algorithm for the k-level uncapacitated facility location problem. Information Processing Letters 72, 161-167

Aardal, K., Labbé, M., Leung, J., Queyranne, M., 1996. On the Two-Level Uncapacitated Facility Location Problem. INFORMS Journal on Computing 8, 289-301

Badejo, O., Ierapetritou, M., 2023. A mathematical modeling approach for supply chain management under disruption and operational uncertainty. AIChE Journal 69, e18037

Baumol, W.J., Wolfe, P., 1958. A Warehouse-Location Problem. Operations Research 6, 252-263

Belieres, S., Hewitt, M., Jozefowiez, N., Semet, F., 2021. A time-expanded network reduction matheuristic for the logistics service network design problem. Transportation Research Part E: Logistics and Transportation Review 147, 102203

Byrka, J., Li, S., Rybicki, B., 2016. Improved Approximation Algorithm for k-level Uncapacitated Facility Location Problem (with Penalties). Theory of Computing Systems 58, 19-44

Cao, Q., Tang, Y., Perera, S., Zhang, J., 2022. Manufacturer- versus retailer-initiated bundling: Implications for the supply chain. Transportation Research Part E: Logistics and Transportation Review 157, 102552

Cortinhal, M.J., Lopes, M.J., Melo, M.T., 2019. A multi-stage supply chain network design problem with in-house production and partial product outsourcing. Applied Mathematical Modelling 70, 572-594

da Silveira Farias, E., Li, J.-Q., Galvez, J.P., Borenstein, D., 2017. Simple heuristic for the strategic supply chain design of large-scale networks: A Brazilian case study. Computers & Industrial Engineering 113, 746-756

Drexl, M.A., 2011. An approximation algorithm for the k-Level Concentrator Location Problem. Operations Research Letters 39, 355-358

Dukkanci, O., Campbell, J.F., Kara, B.Y., 2024. Facility location decisions for drone delivery: A literature review. European Journal of Operational Research 316, 397-418

Ekici, A., Keskinocak, P., Swann, J.L., 2013. Modeling Influenza Pandemic and Planning Food Distribution. Manufacturing & Service Operations Management 16, 11-27

Eskandarpour, M., Dejax, P., Péton, O., 2017. A large neighborhood search heuristic for supply chain network design. Computers & Operations Research 80, 23-37

Farahani, R.Z., Fallah, S., Ruiz, R., Hosseini, S., Asgari, N., 2019. OR models in urban service facility location: A critical review of applications and future developments. European Journal of Operational Research 276, 1-27

Farahani, R.Z., Hekmatfar, M., Fahimnia, B., Kazemzadeh, N., 2014. Hierarchical facility location problem: Models, classifications, techniques, and applications. Computers & Industrial Engineering 68, 104-117

Fattahi, M., Govindan, K., Keyvanshokooh, E., 2017. Responsive and resilient supply chain network design under operational and disruption risks with delivery lead-time sensitive customers. Transportation Research Part E: Logistics and Transportation Review 101, 176-200

Gelman, A., 2005. Analysis of variance—why it is more important than ever. The Annals of Statistics 33, 1-53

Gendron, B., Khuong, P.-V., Semet, F., 2017. Comparison of formulations for the two-level uncapacitated facility location problem with single assignment constraints. Computers & Operations Research 86, 86-93

Gendron, B., Khuong, P.-V., Semet, F., 2023. Models and Methods for Two-Level Uncapacitated Facility Location problems. In: Springer (ed.) Contributions to Combinatorial Optimization and Applications.

Gendron, B., Semet, F., 2009. Formulations and relaxations for a multi-echelon capacitated location–distribution problem. Computers & Operations Research 36, 1335-1355

Guo, Y., Yu, J., Allaoui, H., Choudhary, A., 2022. Lateral collaboration with cost-sharing in sustainable supply chain optimisation: A combinatorial framework. Transportation Research Part E: Logistics and Transportation Review 157, 102593

Hajipour, V., Fattahi, P., Tavana, M., Di Caprio, D., 2016. Multi-objective multi-layer congested facility location-allocation problem optimization with Pareto-based meta-heuristics. Applied Mathematical Modelling 40, 4948-4969

Jahani, H., Abbasi, B., Sheu, J.-B., Klibi, W., 2024. Supply chain network design with financial considerations: A comprehensive review. European Journal of Operational Research 312, 799-839

Janjevic, M., Winkenbach, M., Merchán, D., 2019. Integrating collection-and-delivery points in the strategic design of urban last-mile e-commerce distribution networks. Transportation Research Part E: Logistics and Transportation Review 131, 37-67

Kang, N., Shen, H., Xu, Y., 2021. JD.com Improves Delivery Networks by a Multiperiod Facility Location Model. INFORMS Journal on Applied Analytics 52, 133-148

Kantor, E., Peleg, D., 2009. Approximate hierarchical facility location and applications to the bounded depth Steiner tree and range assignment problems. Journal of Discrete Algorithms 7, 341-362

Karatas, M., Eriskin, L., 2023. Linear and piecewise linear formulations for a hierarchical facility location and sizing problem. Omega 118, 102850

Kaufman, L., Eede, M.V., Hansen, P., 1977. A Plant and Warehouse Location Problem. Journal of the Operational Research Society 28, 547-554

Khalifehzadeh, S., Seifbarghy, M., Naderi, B., 2015. A four-echelon supply chain network design with shortage: Mathematical modeling and solution methods. Journal of Manufacturing Systems 35, 164-175

Klose, A., Drexl, A., 2005. Facility location models for distribution system design. European Journal of Operational Research 162, 4-29

Kratica, J., Dugošija, D., Savić, A., 2014. A new mixed integer linear programming model for the multi level uncapacitated facility location problem. Applied Mathematical Modelling 38, 2118-2129

Kumar, A., Kumar, K., 2024. An uncertain sustainable supply chain network design for regulating greenhouse gas emission and supply chain cost. Cleaner Logistics and Supply Chain 10, 100142

Kumar, R., Ganapathy, L., Gokhale, R., Tiwari, M.K., 2020. Quantitative approaches for the integration of production and distribution planning in the supply chain: a systematic literature review. International Journal of Production Research 58, 3527-3553

Lai, C.-M., Chiu, C.-C., Liu, W.-C., Yeh, W.-C., 2019. A novel nondominated sorting simplified swarm optimization for multi-stage capacitated facility location problems with multiple quantitative and qualitative objectives. Applied Soft Computing 84, 105684

Latha Shankar, B., Basavarajappa, S., Chen, J.C.H., Kadadevaramath, R.S., 2013. Location and allocation decisions for multi-echelon supply chain network – A multi-objective evolutionary approach. Expert Systems with Applications 40, 551-562

Lee, H., Whang, S., 1999. Decentralized Multi-Echelon Supply Chains: Incentives and Information. Management Science 45, 633-640

Li, R., Huang, H.C., Huang, J., 2013. Heuristic algorithms for general k-level facility location problems. Journal of the Operational Research Society 64, 106-113

Majumdar, A., Singh, S.P., Jessica, J., Agarwal, A., 2023. Network design for a decarbonised supply chain considering cap-and-trade policy of carbon emissions. Annals of Operations Research

Malik, A., Contreras, I., Vidyarthi, N., 2022. Two-Level Capacitated Discrete Location with Concave Costs. Transportation Science 56, 1703-1722

Maxwell, S.E., Delaney, H.D., Kelley, K., 2017. Designing experiments and analyzing data: A model comparison perspective. Routledge.

Melo, M.T., Nickel, S., Saldanha-da-Gama, F., 2009. Facility location and supply chain management – A review. European Journal of Operational Research 196, 401-412

Melo, M.T., Nickel, S., Saldanha-da-Gama, F., 2012. A tabu search heuristic for redesigning a multi-echelon supply chain network over a planning horizon. International Journal of Production Economics 136, 218-230

Melo, M.T., Nickel, S., Saldanha-da-Gama, F., 2014. An efficient heuristic approach for a multi-period logistics network redesign problem. TOP 22, 80-108

Melo, M.T., Nickel, S., Saldanha da Gama, F., 2006. Dynamic multi-commodity capacitated facility location: a mathematical modeling framework for strategic supply chain planning. Computers & Operations Research 33, 181-208

Mohammadi, A., Abbasi, A., Alimohammadlou, M., Eghtesadifard, M., Khalifeh, M., 2017. Optimal design of a multi-echelon supply chain in a system thinking framework: An integrated financial-operational approach. Computers & Industrial Engineering 114, 297-315

Mortazavi, A., Arshadi Khamseh, A., Azimi, P., 2015. Designing of an intelligent self-adaptive model for supply chain ordering management system. Engineering Applications of Artificial Intelligence 37, 207-220

Myung, Y.-S., Yu, Y.-M., 2020. Freight transportation network model with bundling option. Transportation Research Part E: Logistics and Transportation Review 133, 101827

Ortiz-Astorquiza, C., Contreras, I., Laporte, G., 2017. Formulations and Approximation Algorithms for Multilevel Uncapacitated Facility Location. INFORMS Journal on Computing 29, 767-779

Ortiz-Astorquiza, C., Contreras, I., Laporte, G., 2018. Multi-level facility location problems. European Journal of Operational Research 267, 791-805

Ortiz-Astorquiza, C., Contreras, I., Laporte, G., 2019. An Exact Algorithm for Multilevel Uncapacitated Facility Location. Transportation Science 53, 1085-1106

Ou-Yang, C., Ansari, R., 2017. Applying a hybrid particle swarm optimization_Tabu search algorithm to a facility location case in Jakarta. Journal of Industrial and Production Engineering 34, 199-212

Ozgen, D., Gulsun, B., 2014. Combining possibilistic linear programming and fuzzy AHP for solving the multi-objective capacitated multi-facility location problem. Information Sciences 268, 185-201

Rafiei, H., Safaei, F., Rabbani, M., 2018. Integrated production-distribution planning problem in a competition-based four-echelon supply chain. Computers & Industrial Engineering 119, 85-99

Ramshani, M., Ostrowski, J., Zhang, K., Li, X., 2019. Two level uncapacitated facility location problem with disruptions. Computers & Industrial Engineering 137, 106089

Ridderstedt, I., Nilsson, J.-E., 2023. Economies of scale versus the costs of bundling: Evidence from procurements of highway pavement replacement. Transportation Research Part A: Policy and Practice 173, 103701

Ruvalcaba-Sandoval, D.A., Olivares-Benitez, E., Rojas, O., Sosa-Gómez, G., 2021. Matheuristics for the Design of a Multi-Step, Multi-Product Supply Chain with Multimodal Transport. In: Applied Sciences

Şahin, G., Süral, H., 2007. A review of hierarchical facility location models. Computers & Operations Research 34, 2310-2331

Saldanha-da-Gama, F., 2022. Facility Location in Logistics and Transportation: An enduring relationship. Transportation Research Part E: Logistics and Transportation Review 166, 102903

Sebatjane, M., Adetunji, O., 2024. A four-echelon supply chain inventory model for growing items with imperfect quality and errors in quality inspection. Annals of Operations Research 335, 327-359

Shoja, A., Molla-Alizadeh-Zavardehi, S., Niroomand, S., 2019. Adaptive meta-heuristic algorithms for flexible supply chain network design problem with different delivery modes. Computers & Industrial Engineering 138, 106107

Sluijk, N., Florio, A.M., Kinable, J., Dellaert, N., Van Woensel, T., 2023. Two-echelon vehicle routing problems: A literature review. European Journal of Operational Research 304, 865-886

Soshko, O., Merkuryev, Y., Chakste, M., 2007. Application in Retail: Locating a Distribution Center. In: Chandra C & Grabis J (eds.) Supply Chain Configuration: Concepts, Solutions, and Applications. Springer US, Boston, MA, pp. 303-333.

Tang, Y., Perera, S.C., Cao, Q., Ji, X., 2023. Supplier versus platform bundling: Optimal strategies under agency selling. Transportation Research Part E: Logistics and Transportation Review 179, 103325

Vanovermeire, C., Sörensen, K., 2014. Integration of the cost allocation in the optimization of collaborative bundling. Transportation Research Part E: Logistics and Transportation Review 72, 125-143

Vishnu, C.R., Das, S.P., Sridharan, R., Ram Kumar, P.N., Narahari, N.S., 2021. Development of a reliable and flexible supply chain network design model: a genetic algorithm based approach. International Journal of Production Research 59, 6185-6209

Wang, F.-M., Wang, J.-J., Li, N., Jiang, Y.-J., Li, S.-C., 2022. A Cost-Sharing Scheme for the k-Level Facility Location Game with Penalties. Journal of the Operations Research Society of China 10, 173-182

Wu, T., Xiao, F., Zhang, C., Zhang, D., Liang, Z., 2019. Regression and extrapolation guided optimization for production–distribution with ship–buy–exchange options. Transportation Research Part E: Logistics and Transportation Review 129, 15-37

Yan, Y., Chow, A.H.F., Ho, C.P., Kuo, Y.-H., Wu, Q., Ying, C., 2022. Reinforcement learning for logistics and supply chain management: Methodologies, state of the art, and future opportunities. Transportation Research Part E: Logistics and Transportation Review 162, 102712

Zhang, J., Ye, Y., 2002. A note on the maximization version of the multi-level facility location problem. Operations Research Letters 30, 333-335

Zokaee, S., Jabbarzadeh, A., Fahimnia, B., Sadjadi, S.J., 2017. Robust supply chain network design: an optimization model with real world application. Annals of Operations Research 257, 15-44

Tang, Y., Perera, S.C., Cao, Q., Ji, X., 2023. Supplier versus platform bundling: Optimal strategies under agency selling. Transportation Research Part E., 179, 103325

Tapia-Ubeda, F.J, Miranda-Gonzalez, P.A, and Gutierrez-Jarpa, G. 2024. Integrating supplier selection decisions into an inventory location problem for designing the supply chain network., Journal of Combinatorial Optimization (2024) 47:2.

Vanovermeire, C., Sörensen, K., 2014. Integration of the cost allocation in the optimization of collaborative bundling. Transportation Research Part E., 72, 125-143

Wang, H, and Alidaee, B., 2019, Effective heuristic for large-scale unrelated parallel machines scheduling problems., Omega, 83, 261-274.

Wang, H, and Alidaee, B., 2023. A new hybrid-heuristic for large-scale combinatorial optimization: A case of quadratic assignment problem., Computers & Industrial Engineering, 179, 109220

Wang, K-J, and Febri, N. 2024. The vending machine deployment and shelf display problem: A bi-layer optimization approach., Transportation Research Part E., 187, 103581.

Wang, D, Yang, K, Yuen, K.F., Yang, L., and Dong, J. 2024. Hybrid risk-averse location-inventory-allocation with secondary disaster considerations in disaster relief logistics: A distributionally robust approach., Transportation Research Part E., 186, 103558.

Wang, F.-M., Wang, J.-J., Li, N., Jiang, Y.-J., Li, S.-C., 2022. A Cost-Sharing Scheme for the k-Level Facility Location Game with Penalties. Journal of the Operations Research Society of China 10, 173-182

Wu, T., Xiao, F., Zhang, C., Zhang, D., Liang, Z., 2019. Regression and extrapolation guided optimization for production–distribution with ship–buy–exchange options. Transportation Research Part E., 129, 15-37

Yan, Y., Chow, A.H.F., Ho, C.P., Kuo, Y.-H., Wu, Q., Ying, C., 2022. Reinforcement learning for logistics and supply chain management: Methodologies, state of the art, and future opportunities. Transportation Research Part E., 162, 102712

Zarandi, R.Z, Motlaq-Kashani, A.S, and Sheikhalishahi, M. 2024. A two-echelon sustainable multi-route location routing problem for biomass supply chain network design considering disruption., Computers and Chemical Engineering, 187, 108744.

Zhang, J., Ye, Y., 2002. A note on the maximization version of the multi-level facility location problem. Operations Research Letters 30, 333-335

Zokaee, S., Jabbarzadeh, A., Fahimnia, B., Sadjadi, S.J., 2017. Robust supply chain network design: an optimization model with real world application. Annals of Operations Research 257, 15-44

Table 1: Results of 15 runs for Local Search with (LS W/Seq) and without (No/Seq) sequence diversification implemented and different size problems

| Problem ID | BFS | | TB (Seconds) | |
|---|---|---|---|---|
| | LS No/Seq | LS W/Seq | LS No/Seq | LS W/Seq |
| MFL-2000-30-50-150-1 | 61630 | 62122 | 1726.348 | 181.595 |
| MFL-2000-30-50-150-2 | 62492 | 63040 | 3745.107 | 73.286 |
| MFL-2000-30-50-150-3 | 61507 | 61646 | 1852.425 | 221.44 |
| MFL-2000-30-50-150-4 | 61245 | 61921 | 3830.775 | 44.141 |
| MFL-2000-30-50-150-5 | 61930 | 62176 | 3527.513 | 315.984 |
| MFL-2000-30-50-150-6 | 61782 | 61927 | 3229.965 | 160.065 |
| MFL-2000-30-50-150-7 | 63199 | 63633 | 2057.345 | 421.592 |
| MFL-2000-30-50-150-8 | 61697 | 62472 | 1280.376 | 345.35 |
| MFL-2000-30-50-150-9 | 62101 | 62734 | 2708.028 | 65.963 |
| MFL-2000-30-50-150-10 | 61361 | 61803 | 3331.392 | 435.053 |
| MFL-3000-30-50-150-1 | 92490 | 93257 | 2523.756 | 682.149 |
| MFL-3000-30-50-150-2 | 94259 | 94726 | 1470.483 | 823.921 |
| MFL-3000-30-50-150-3 | 93604 | 94611 | 4769.702 | 414.944 |
| MFL-3000-30-50-150-4 | 94326 | 95599 | 4507.413 | 207.148 |
| MFL-3000-30-50-150-5 | 92679 | 93344 | 1519.84 | 117.414 |
| MFL-3000-30-50-150-6 | 93463 | 93880 | 327.474 | 694.27 |
| MFL-3000-30-50-150-7 | 93270 | 93473 | 2084.077 | 513.646 |
| MFL-3000-30-50-150-8 | 93946 | 94700 | 2725.755 | 166.2 |
| MFL-3000-30-50-150-9 | 93519 | 94572 | 7797.316 | 849.297 |
| MFL-3000-30-50-150-10 | 92776 | 93736 | 5954.239 | 490.718 |
| MFL-4000-30-50-150-1 | 125858 | 126888 | 15123.270 | 443.797 |
| MFL-4000-30-50-150-2 | 124319 | 126047 | 6096.143 | 957.803 |
| MFL-4000-30-50-150-3 | 124771 | 125631 | 7479.600 | 1354.779 |
| MFL-4000-30-50-150-4 | 126006 | 126211 | 1999.664 | 804.443 |
| MFL-4000-30-50-150-5 | 124771 | 124885 | 2908.508 | 605.723 |
| MFL-4000-30-50-150-6 | 126014 | 126755 | 9519.062 | 1097.527 |
| MFL-4000-30-50-150-7 | 125199 | 125990 | 11586.715 | 1160.502 |
| MFL-4000-30-50-150-8 | 124862 | 125248 | 4058.730 | 961.35 |
| MFL-4000-30-50-150-9 | 124012 | 125406 | 8898.456 | 1910.242 |
| MFL-4000-30-50-150-10 | 126551 | 127766 | 14458.103 | 646.801 |
| MFL-5000-30-50-150-1 | 157136 | 158286 | 3804.742 | 2629.553 |
| MFL-5000-30-50-150-2 | 155019 | 155735 | 15103.743 | 1112.674 |
| MFL-5000-30-50-150-3 | 154729 | 156398 | 21890.742 | 2855.883 |
| MFL-5000-30-50-150-4 | 156106 | 156944 | 2419.119 | 2272.559 |
| MFL-5000-30-50-150-5 | 157244 | 159051 | 2580.811 | 1142.879 |
| MFL-5000-30-50-150-6 | 155479 | 155674 | 18343.008 | 798.002 |
| MFL-5000-30-50-150-7 | 156677 | 158414 | 24297.361 | 1275.738 |
| MFL-5000-30-50-150-8 | 154841 | 156788 | 18190.207 | 1836.021 |
| MFL-5000-30-50-150-9 | 156403 | 156583 | 24919.613 | 2352.625 |
| MFL-5000-30-50-150-10 | 156635 | 158449 | 5207.167 | 3011.262 |
| MFL-6000-30-50-150-1 | 189058 | 190256 | 10927.647 | 2488.191 |
| MFL-6000-30-50-150-2 | 188377 | 188879 | 32426.342 | 1257.266 |
| MFL-6000-30-50-150-3 | 187405 | 188998 | 28971.986 | 4261.887 |
| MFL-6000-30-50-150-4 | 190525 | 191049 | 23910.098 | 838.773 |
| MFL-6000-30-50-150-5 | 188544 | 190849 | 32917.391 | 3459.965 |
| MFL-6000-30-50-150-6 | 185374 | 188390 | 33451.566 | 972.449 |
| MFL-6000-30-50-150-7 | 189359 | 190430 | 5368.320 | 4345.902 |
| MFL-6000-30-50-150-8 | 189956 | 190555 | 28552.289 | 437.605 |
| MFL-6000-30-50-150-9 | 189688 | 191025 | 21427.699 | 3856.262 |
| MFL-6000-30-50-150-10 | 189667 | 191850 | 16234.212 | 1713.406 |

Table 2: Results for Local Search and Tabu Search, both with sequence diversification implemented and different size problems

| Problem ID | BFS | | TB (Seconds) | |
|---|---|---|---|---|
| | LS W/Seq | TS W/Seq | LS W/Seq | TS W/Seq |
| MFL-2000-30-50-150-1 | 61629 | 62329 | 131.25 | 93.107 |
| MFL-2000-30-50-150-2 | 62841 | 63357 | 320.66 | 24.316 |
| MFL-2000-30-50-150-3 | 62537 | 64201 | 384.787 | 119.711 |
| MFL-2000-30-50-150-4 | 62390 | 63006 | 326.883 | 83.815 |
| MFL-2000-30-50-150-5 | 62114 | 63113 | 329.955 | 48.839 |
| MFL-2000-30-50-150-6 | 63911 | 65038 | 406.71 | 15.756 |
| MFL-2000-30-50-150-7 | 61875 | 62493 | 119.463 | 36.981 |
| MFL-2000-30-50-150-8 | 61723 | 62476 | 227.707 | 45.674 |
| MFL-2000-30-50-150-9 | 61675 | 63775 | 74.207 | 34.556 |
| MFL-2000-30-50-150-10 | 61884 | 62531 | 428.186 | 73.378 |
| MFL-3000-30-50-150-1 | 95041 | 96131 | 379.042 | 97.79 |
| MFL-3000-30-50-150-2 | 93043 | 94122 | 572.949 | 109.251 |
| MFL-3000-30-50-150-3 | 94205 | 95107 | 115.479 | 62.744 |
| MFL-3000-30-50-150-4 | 93262 | 95764 | 647.689 | 90.668 |
| MFL-3000-30-50-150-5 | 94185 | 94914 | 299.48 | 36.014 |
| MFL-3000-30-50-150-6 | 93333 | 94209 | 146.055 | 84.956 |
| MFL-3000-30-50-150-7 | 95676 | 96394 | 120.55 | 60.992 |
| MFL-3000-30-50-150-8 | 93535 | 94590 | 534.518 | 124.549 |
| MFL-3000-30-50-150-9 | 94285 | 95525 | 605.625 | 93.656 |
| MFL-3000-30-50-150-10 | 94667 | 95420 | 386.158 | 66.653 |
| MFL-4000-30-50-150-1 | 125664 | 130818 | 570.113 | 696.824 |
| MFL-4000-30-50-150-2 | 124784 | 126063 | 2068.355 | 201.914 |
| MFL-4000-30-50-150-3 | 128473 | 132390 | 275.664 | 192.693 |
| MFL-4000-30-50-150-4 | 125349 | 126661 | 1509.285 | 169.121 |
| MFL-4000-30-50-150-5 | 125503 | 126763 | 1502.961 | 113.737 |
| MFL-4000-30-50-150-6 | 125517 | 129593 | 1453.544 | 204.998 |
| MFL-4000-30-50-150-7 | 125740 | 126993 | 481.353 | 145.972 |
| MFL-4000-30-50-150-8 | 126055 | 127458 | 896.573 | 279.699 |
| MFL-4000-30-50-150-9 | 126233 | 130044 | 414.342 | 393.488 |
| MFL-4000-30-50-150-10 | 125837 | 127767 | 439.529 | 154.405 |
| MFL-5000-30-50-150-1 | 156681 | 160816 | 614.694 | 273.674 |
| MFL-5000-30-50-150-2 | 157187 | 165447 | 1274.16 | 227.326 |
| MFL-5000-30-50-150-3 | 158747 | 161427 | 240.041 | 218.206 |
| MFL-5000-30-50-150-4 | 157219 | 158720 | 1815.309 | 315.53 |
| MFL-5000-30-50-150-5 | 158110 | 159451 | 2444.637 | 267.499 |
| MFL-5000-30-50-150-6 | 157757 | 159514 | 753.526 | 377.615 |
| MFL-5000-30-50-150-7 | 158237 | 159606 | 1188.45 | 159.828 |
| MFL-5000-30-50-150-8 | 157732 | 160910 | 1996.789 | 298.26 |
| MFL-5000-30-50-150-9 | 156957 | 159475 | 1024.633 | 245.323 |
| MFL-5000-30-50-150-10 | 156929 | 159850 | 2025.307 | 321.59 |
| MFL-6000-30-50-150-1 | 188087 | 191476 | 2492.757 | 270.323 |
| MFL-6000-30-50-150-2 | 189640 | 191381 | 596.487 | 311.141 |
| MFL-6000-30-50-150-3 | 188109 | 190161 | 1501.423 | 433.402 |
| MFL-6000-30-50-150-4 | 190034 | 191799 | 879.954 | 442.25 |
| MFL-6000-30-50-150-5 | 186807 | 188731 | 1006.758 | 341.676 |
| MFL-6000-30-50-150-6 | 187198 | 188521 | 1633.959 | 195.176 |
| MFL-6000-30-50-150-7 | 190927 | 192305 | 3092.023 | 1432.937 |
| MFL-6000-30-50-150-8 | 189214 | 192464 | 1814.912 | 420.329 |
| MFL-6000-30-50-150-9 | 191293 | 193905 | 2843.013 | 92.682 |
| MFL-6000-30-50-150-10 | 191798 | 193021 | 2296.57 | 62.328 |

Table 3 Descriptive statistics of LS/TS for consistency evaluation across 15 runs on all instances

| Local Search with (LS W/Seq) and without (No/Seq) sequence | | | | | | |
|---|---|---|---|---|---|---|
| Problem | Number of Instances | Runs per instance | BFS_diff | | TB_diff | |
| | | | Mean | STDEV | Mean | STDEV |
| 2000-30-50-150 | 10 | 15 | 453 (0.73%) | 208.12 | -2502.48 (92%) | 943.43 |
| 3000-30-50-150 | 10 | 15 | 757 (0.81%) | 311.89 | -2872.03 (85%) | 2191.43 |
| 4000-30-50-150 | 10 | 15 | 846 (0.68%) | 492.02 | -7218.53 (88%) | 4403.09 |
| 5000-30-50-150 | 10 | 15 | 1205 (0.77%) | 650.36 | -11746.93 (86%) | 8949.94 |
| 6000-30-50-150 | 10 | 15 | 1433 (0.76%) | 802.40 | -21055.58 (90%) | 9929.54 |
| **Tabu Search and Local Search, both with sequence** | | | | | | |
| 2000-30-50-150 | 10 | 15 | 974 (1.6%) | 496.34 | -217.37 (79%) | 120.58 |
| 3000-30-50-150 | 10 | 15 | 1094 (1.8%) | 498.08 | -298.03 (78%) | 181.09 |
| 4000-30-50-150 | 10 | 15 | 2540 (2.0%) | 1440.89 | -705.89 (73%) | 662.07 |
| 5000-30-50-150 | 10 | 15 | 2966 (1.9%) | 1962.69 | -1067.27 (80%) | 663.75 |
| 6000-30-50-150 | 10 | 15 | 2066 (1.1%) | 735.05 | -1415.56 (78%) | 779.79 |

Note. For Local Search with and without sequence, *BFS_diff* represents the difference between the BFS value of Local Search with sequence and that of Local Search without sequence, while *TB_diff* is calculated similarly for the TB value. For the comparison between Tabu Search and Local Search, *BFS_diff* denotes the difference between the BFS value of Tabu Search with sequence and that of Local Search with sequence, and *TB_diff* is computed similarly.

Table 4. ANOVA and t-test results of algorithms with and without sequence

| Local Search with (LS W/Seq) and without (No/Seq) sequence | | | | | | |
|---|---|---|---|---|---|---|
| Problem | Number of instances | Runs per instance | BFS | | TB | |
| | | | ANOVA F-value | T-test t-value | ANOVA F-value | T-test t-value |
| 2000-30-50-150 | 10 | 15 | 2.80. | -6.53*** | 71.07*** | 7.96*** |
| 3000-30-50-150 | 10 | 15 | 5.87* | -7.28*** | 13.26*** | 3.93*** |
| 4000-30-50-150 | 10 | 15 | 4.15. | -5.16*** | 21.37*** | 4.92*** |
| 5000-30-50-150 | 10 | 15 | 5.61* | -5.56*** | 18.03*** | 3.94*** |
| 6000-30-50-150 | 10 | 15 | 4.90* | -5.36*** | 50.85*** | 6.36*** |
| **Tabu Search and Local Search, both with sequence** | | | | | | |
| 2000-30-50-150 | 10 | 15 | 7.16** | -5.89*** | 30.30*** | 5.41*** |
| 3000-30-50-150 | 10 | 15 | 8.66** | -6.59*** | 16.68*** | 4.94*** |
| 4000-30-50-150 | 10 | 15 | 8.53** | -5.29*** | 13.59*** | 3.20*** |
| 5000-30-50-150 | 10 | 15 | 16.00*** | -4.53*** | 23.68*** | 4.82*** |
| 6000-30-50-150 | 10 | 15 | 4.97* | -8.43*** | 17.24*** | 5.45*** |

Note: (. p-value<=0.1, * p-value<=0.05, ** p-value<=0.01, ***, p-value<=0.005)

Table 5. **Wilcoxon signed-rank test (p-value)** results of algorithms with and without sequence

| **Local Search with (LS W/Seq) and without (No/Seq) sequence** | | | | |
|---|---|---|---|---|
| Problem | Number of instances | Runs per instance | BFS | TB |
| 2000-30-50-150 | 10 | 15 | 0.038* | 5.41E-06*** |
| 3000-30-50-150 | 10 | 15 | 0.026* | 2.44E-04*** |
| 4000-30-50-150 | 10 | 15 | 0.018* | 5.41E-06*** |
| 5000-30-50-150 | 10 | 15 | 0.032* | 1.62E-04*** |
| 6000-30-50-150 | 10 | 15 | 0.014* | 5.41E-06*** |
| **Tabu Search and Local Search, both with sequence** | | | | |
| 2000-30-50-150 | 10 | 15 | 5.75E-03** | 6.50E-05*** |
| 3000-30-50-150 | 10 | 15 | 7.34E-03** | 2.17E-05*** |
| 4000-30-50-150 | 10 | 15 | 2.44E-04*** | 2.44E-04*** |
| 5000-30-50-150 | 10 | 15 | 1.08E-05*** | 2.44E-04*** |
| 6000-30-50-150 | 10 | 15 | 5.75E-03** | 3.79E-05*** |

Note: (. p-value<=0.1, * p-value<=0.05, ** p-value<=0.01, ***, p-value<=0.005)

## APPENDIX A

The following proposition (*n*, *k*, *j*, *m*) means if retail store *m* receives its supply from plant *n*, which is served by a warehouse *k*, and a distribution center *i*.

**Proposition 1:** Given a solution (*x*, *y*), considering the following simple exchange moves, the amount of change in the objective function in each case is calculated as follows:

(a) An open retail store *m* currently served via the path (*n*, *k*, *j*, *m*) is closed. The amount of change is equal to Change(a).

$$\text{Change}(a) = -R_m + c_{jm} + c_{kj} + c_{nk} + fs_m + \Delta_1 + \Delta_2 + \Delta_3$$

$\Delta_1 = fd_j$, if retail store *m* is the only store that the distribution center *j* is currently serving, 0 otherwise

$\Delta_2 = fw_k$ , if *m* is the only store that the warehouse *k* is currently serving, 0 otherwise

$\Delta_3 = fp_n$, if *m* is the only store that the plant *n* is currently serving, 0 otherwise

(b) A closed retail store *m* is opened and served via a path (*n*, *k*, *j*, *m*). The amount of change is equal to Change(b).

$$\text{Change(b)} = R_m - c_{jm} - c_{kj} - c_{nk} - fs_m + \Delta_1 + \Delta_2 + \Delta_3$$

$\Delta_1 = -fd_j$ , if *m* is the only store that the distribution center *j* is currently serving, 0 otherwise

$\Delta_2 = -fw_k$ , if *m* is the only store that the warehouse *k* is currently serving, 0 otherwise

$\Delta_3 = -fp_n$ , if *m* is the only store that the plant *n* is currently serving, 0 otherwise

(c) A retail store *m* is served via a distribution center *j'* instead of the distribution center *j*, **Figure 2(a)**, the amount of change is equal to Change (c).

$$\text{Change(c)} = -c_{jm} + c_{j'm} - c_{kj} + c_{kj'} + \Delta_1 + \Delta_2$$

$\Delta_1 = -fd_j$ , if *m* is the only store that the distribution center *j* is currently serving, 0 otherwise

$\Delta_2 = fd_{j'}$, if *m* is the only store that the distribution center *j'* will be serving, 0 otherwise

(d) A retail store *m* is served via a warehouse *k'* instead of warehouse *k*, **Figure2(b)**, the amount of change is equal to Change (d).

$$\text{Change(d)} = -c_{nk} + c_{nk'} - c_{kj} + c_{k'j} + \Delta_1 + \Delta_2$$

$\Delta_1 = -fw_k$ , if *m* is the only store that the warehouse *k* is currently serving, 0 otherwise

$\Delta_2 = fw_{k'}$, if *m* is the only store that the warehouse *k'* will be serving, 0 otherwise

(e) A retail store *m* is served via plant *n'* instead of plant *n*, **Figure2(c)**, the amount of change is equal to Change (e).

$$\text{Change(e)} = -c_{nk} + c_{n'k} + \Delta_1 + \Delta_2$$

$\Delta_1 = -fp_n$ , if *m* is the only store that the plant *n* is currently serving, 0 otherwise

$\Delta_2 = fp_{n'}$, if *m* is the only store that the plant *n'* will be serving, 0 otherwise

**Proof**. We establish the claim by considering the fixed costs associated with each facility and the variable costs of product movement along each arc. Firstly, each plant, warehouse, distribution center, and retail store incur a one-time fixed cost if and only if it is opened. Moreover, the cost of transporting a bundle of products along each arc is store-specific. Therefore, if a retail store m is currently served via the path (*n*, *k*, *j*, *m*), any modification to the plants, warehouses, or distribution centers involved in this path will impact the cost of product movement along the affected arcs. Similarly, if a retail store is closed, all costs associated with the arcs serving that store will be eliminated. Additionally, the revenue generated by serving that store will also be lost.

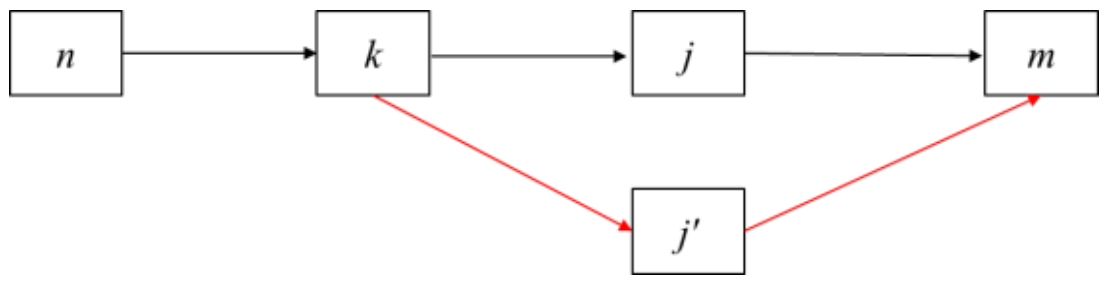

Figure 2(a). Change in the objective function if store *m* is served via path (*n*, *k*, *j'*, *m*) instead of (*n*, *k*, *j*, *m*).

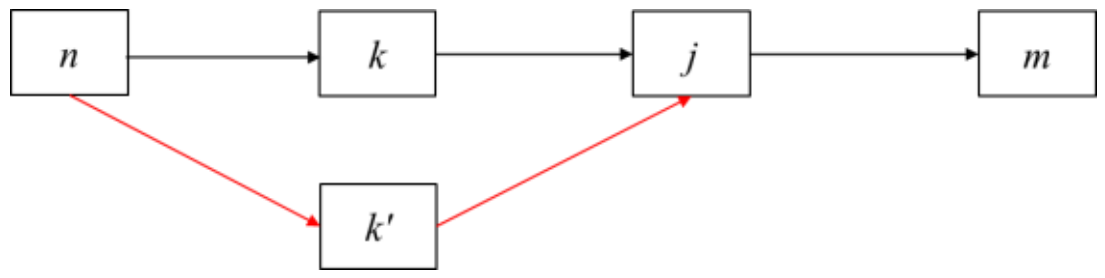

Figure 2(b). Change in the objective function if store *m* is served via path (*n*, *k'*, *j*, *m*) instead of (*n*, *k*, *j*, *m*).

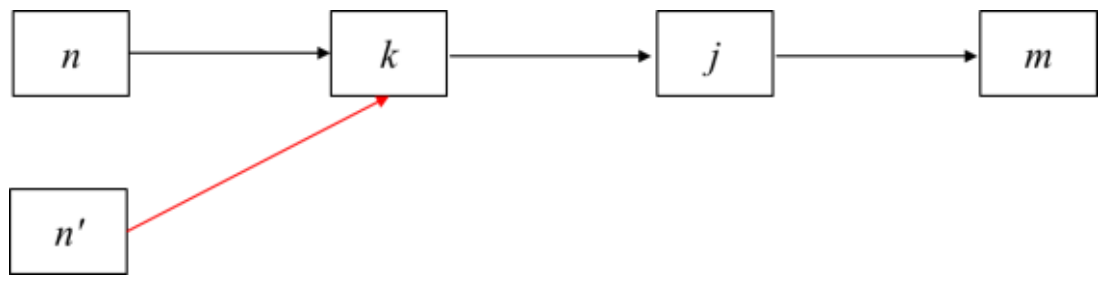

Figure 2(c). Change in the objective function if store *m* is served via path (*n'*, *k*, *j*, *m*) instead of (*n*, *k*, *j*, *m*).

## Appendix B

**Heuristic 1: Simple Local Search with *r-Opt* Sequencing Diversification Strategy**

**Initialization:** Sets: *P*, *W*, *D*, *S*. Vector of fixed costs: *fs*, *fp*, *fw*, *fd*, Upper bounds: *us*, *up*, *uw*, *ud*. Matrix of transportation costs (arc costs) between (*P* and *W*), (*W* and *D*), (*D* and *S*).

Random sequences *LS*, *LP*, *LW*, *LD*, of elements of *S*, *P*, *W*, *and D*, respectively. A feasible solution *is z= (x, y, w, d);* here, there is a path for each opened store m (*n*, *k*, *j*, *m*). Let *z*=z*, and *(n*, k*, j*, m) = (n, k, j, m)* for all opened *m* be the best solution found so far, and keep track of the best solution found throughout the process.

1. Flag=1
2. DO WHILE(Flag=1)
3. Flag=0

4. FOR *m=1*, |*S*|
5. *g=LS(m)*
6. IF (*g* is open) THEN
7. IF (Closing *g* improves the objective function, Case (a), Proposition 1) THEN
8. **Update:** Close *g, z*, (n*, k*, j*, m)*, Flag=1, Go to 25
9. END IF
10. ELSE
11. FOR *j=1*, |*D*|
12. *h=LD(j)*
13. FOR *k=1*, |*W*|
14. *e=LW(k)*
15. FOR *n=1*, |*P*|
16. *f=LP(n)*
17. IF (Opening store *g* is feasible, i.e., (*f, e, h, g*) is feasible, and improves the objective function, Case (b), Proposition 1) THEN
18. **Update:** Open *g, z*, (n*, k*, j*, m)*, Flag=1, Go to 25
19. END IF
20. END FOR
21. END FOR
22. END FOR
23. END IF
24. END FOR
25. Call *r-Opt(.)*
26. FOR *m=1*, |*S*|
27. *g=LS(m)*
28. IF (*g* is open) THEN
29. FOR *j=1*, |*D*|
30. *h=LD(j)*
31. IF ( $y_{hg} = 1$ and $y_{h*g} = 0$ ) is feasible and improves the objective function, Case (c), Proposition 1) THEN
32. **Update:** $y_{hg} = 1$, $y_{h*g} = 0$, *z*, (n*, k*, j*, m)*, Flag=1, Go to 36
33. END IF
34. END FOR

35. END IF
36. IF (*g* is open) THEN
37. FOR *k=1*, |*W*|
38. *e=LW(k)*
39. IF ( $x_{gn^*e} = 1$ and $x_{gn^*k^*} = 0$ ) is feasible and improves the objective function, Case (d), Proposition 1) THEN
40. **Update:** $x_{gn^*e} = 1$, $x_{gn^*k^*} = 0$, *z*, (n*, k*, j*, m),* Flag=1, Go to 44
41. END IF
42. END FOR
43. END IF
44. IF (*g* is open) THEN
45. FOR *n=1*, |*P*|
46. *f=LP(n)*
47. IF ( $x_{gfk^*} = 1$ and $x_{gn^*k^*} = 0$ ) is feasible and improves the objective function, Case (e), Proposition 1) THEN
48. **Update:** $x_{gfk^*} = 1$, $x_{gn^*k^*} = 0$, *z*, (n*, k*, j*, m),* Flag=1, Go to 53
49. END IF
50. END FOR
51. END IF
52. END FOR
53. Call *r-Opt(.)*
54. END WHILE

In Heuristic 1, specifically in Steps 25 and 53, new sequences are generated at each iteration. This can be achieved through various methods. Our approach adopts an r-Opt strategy, commonly employed in heuristics for solving sequencing problems. We implemented a limited double bridge for *r*=4, as introduced in (Alidaee & Wang, 2017). By utilizing different sequences, we can explore more diverse regions of the solution space, effectively implementing a multi-start strategy. Notably, while Cases (a-e) of Proposition 1 are presented in a specific order for implementation simplicity, they can be executed in varying orders, which may lead to even more effective outcomes in practice.

**Heuristic 2: A multi-start Tabu Search with *r-Opt* Sequencing Diversification Strategy**

**Initialization:** Sets: *P, W, D, S*. Vector of fixed costs: *fs, fp, fw, fd,* Upper bounds: *us, up, uw, ud.* Matrix of transportation costs (arc costs) between (*P* and *W*), (*W* and *D*), (*D* and *S*). A criterion to stop, e.g., *MAXCOUNT*.

1. DO COUNT=1, MAXCOUNT
2. Find a new starting feasible solution along the use of sequences found in Step 5:
   *z1=(x1,y1,w1,d1)*, path *(n1,k1,j1,m)* associated with *z1* for all $m \in S$ ,
   calculate the value of objective function *TC1*.
3. Call *Heuristic 3(.)* (This is similar to Heuristic 1 with Tabu strategy implemented)
4. **Update:** *z_best=z*, TC_best=TC*, (n,k,j,m)_best=(n*,k*,j*,m)*
5. CALL *r-Opt(.),*
   (This randomly chooses new sequences *LS, LP, LW, LD*)
7. END DO

Heuristic 3 combines a tabu search with an embedded *r-Opt* sequence diversification strategy. While the tabu strategy could be applied to plants, warehouses, distribution centers, and retail stores, we limited its implementation to retail stores. The reason for this decision is that the fixed costs associated with opening facilities in the other layers (plants, warehouses, and distribution centers) cause the tabu strategy to be ineffective when applied to these facilities. This is because each of these facilities may be utilized multiple times, whereas a retail store is either used or not, making the tabu strategy more suitable for the retail store layer. Note that Heuristic 3 is utilized in the third step of Heuristic 2.

**Heuristic 3: Tabu Search with *r-Opt* Sequencing Diversification Strategy**

**Initialization:** Sets: *P, W, D, S*. Vector of fixed costs: *fs, fp, fw, fd,* Upper bounds: *us, up, uw, ud.* Matrix of transportation costs (arc costs) between (*P* and *W*), (*W* and *D*), (*D* and *S*). Let *z*=z1=(x*,y*,w*,d*)*, and associated path *(n*,k*,j*,m)=(n1,k1,j1,m)* for all *m* and *TC*=TC1*, be the best solution found so far. Random sequences *LP*, *LW*, *LD*, and *LS* of elements *P*, *W*, *D*, and *S*, respectively. Tabu vector: *Tabu_S*. *Tabu_tenure*, and a criterion to stop, e.g., *MAXCOUNT*.

However, Heuristic 3 shares similarities with Heuristic 1, with the key distinction being incorporating a Tabu Strategy specifically designed for retail stores. To transform Heuristic 1 into Heuristic 3, the following modifications are necessary. Additionally, it is essential to update the step numbering accordingly.

**Steps 6-9 of Heuristic 1 should be changed to:**

IF (*g* is open) THEN

Calculate associated *TC* if *g* is closed

IF ((*TC*>*TC1*).and.*Tabu_S(g)=0*).or.(*TC*>*TC**)) THEN

**Update:** Close *g, z1, TC1, (n1, k1, j1, m), z*, (n*, k*, j*, m)*, *Tabu_S*, Flag=1,

**Go to the next appropriate Step**

END IF

**Steps 17-19 of Heuristic 1 should be changed to:**

IF(*g* is closed) THEN

IF(Opening *g* is feasible, i.e., (*f*,*e*,*h*,*g*) is feasible) THEN

Calculate associated *TC* if *g* is opened

IF((*TC*>*TC1*).and.*Tabu_S(g)=0*).or.(*TC*>*TC**)) THEN

**Update:** Open *g, z1,TC1,(n1,k1,j1,m), z*, (n*,k*,j*,m), Tabu_S,* Flag=1,

**Go to the next appropriate Step**

END IF

END IF